\documentclass[12pt]{article}
\title{Geodesic completeness for some meromorphic metrics}
\author{Claudio Meneghini}
\begin{document}
\bibliographystyle{plain}
\maketitle
\def\quan{\vrule height6pt width6pt depth0pt}
\def\QUAN{\nobreak $\ $\quan\hfill\vskip0.6truecm\par}
\def\BETA{\mathop{\beta}\limits}
\def\GAMMA{\mathop{\gamma}\limits}
\def\VI{\mathop{v}\limits}
\def\UI{\mathop{u}\limits}
\def\VII{\mathop{V}\limits}
\def\WI{\mathop{w}\limits}
\def\ZETA{\mathop{Z}\limits}
\def\ssqrt#1{\left(#1\right)^{1/2}}
\def\sssqrt#1{\left(#1\right)^{-1/2}}
\def\TTT{\sl}
\def\BBB{\sl}

\newtheorem{definition}{Definition}[section]
\newtheorem{lemma}[definition]{Lemma}
\newtheorem{proposition}[definition]{Proposition}
\newtheorem{theorem}[definition]{Theorem}        
\newtheorem{corollary}[definition]{Corollary}  
\newtheorem{remark}[definition]{Remark}  

\font\sdopp=msbm10
\def\DI{\sdopp{\hbox{D}}}
\def\ESSE {\sdopp {\hbox{S}}}
\def\ERRE {\sdopp {\hbox{R}}}
\def\CI {\sdopp {\hbox{C}}}
\def\ENNE{\sdopp {\hbox{N}}}
\def\ZETA{\sdopp {\hbox{Z}}}
\def\PI {\sdopp {\hbox{P}}}
\def\M{\hbox{\boldmath{}$M$\unboldmath}} 
\def\N{\hbox{\boldmath{}$N$\unboldmath}} 
\def\P{\hbox{\boldmath{}$P$\unboldmath}} 
\def\Y{\hbox{\boldmath{}$Y$\unboldmath}} 
\def\tr{\hbox{\boldmath{}$tr$\unboldmath}} 
\def\f{\hbox{\boldmath{}$f$\unboldmath}} 
\def\u{\hbox{\boldmath{}$u$\unboldmath}} 
\def\v{\hbox{\boldmath{}$v$\unboldmath}} 
\def\U{\hbox{\boldmath{}$U$\unboldmath}} 
\def\V{\hbox{\boldmath{}$V$\unboldmath}} 
\def\W{\hbox{\boldmath{}$W$\unboldmath}} 
\def\id{\hbox{\boldmath{}$id$\unboldmath}} 
\def\alph{\hbox{\boldmath{}$\aleph$\unboldmath}} 
\def\bet{\hbox{\boldmath{}$\beta$\unboldmath}} 
\def\gam{\hbox{\boldmath{}$\gamma$\unboldmath}} 
\def\U{\mathop{u}\limits}
\def\f{\hbox{\boldmath{}$f$\unboldmath}} 
\def\g{\hbox{\boldmath{}$g$\unboldmath}} 
\def\h{\hbox{\boldmath{}$h$\unboldmath}} 
\def\IM{\hbox{\boldmath{}$i$\unboldmath}} 

\sloppy
\parindent=0pt
\section{Foreword}
In this paper we shall be concerned with generalizing the ideas of 'metric' and geodesic
for a complex manifold $\M$: we emphasize
that our curves will be complex ones; a metric
will be, informally speaking, a {\it symmetric}
quadratic form on the holomorphic tangent space
at each point $p\in\M$, holomorphically
depending on the point itself; of course, it couldn't
have any 'signature', but, by simmetry, it induces
a canonical Levi-Civita's connexion on $\M$,
which in turn allows us to define geodesics to
be auto-parallel paths.
We illustrate some motivations  
(see \cite{dubnovfom} p.186 ff):
consider the space 
${\cal F}$ of antisymmetric covariant tensors of rank two
in Minkowski's space $\ERRE_{1,3}$: 
electromagnetic fields are such ones.
Let $F\in{\cal F}$: we can write $F=
\sum_{i<j}F_{ij}dx^i\wedge dx^j$ where 
$x^0...x^3$ are the natural coordinate functions
on $\ERRE_{1,3}$.
At each point, the space ${\cal F}_p$ of 
all tensors in ${\cal F}$ evaluated at $p$ is a six-
dimensional real vector space; moreover, the
adjoint operator $*$ with respect to Minkowski's metric is such that $**=-1$:
all these facts imply that 
${\cal F}_p$ could be thought of as a complex three dimensional
vector space ${\cal G}_p$ by setting $(a+\IM b)F=aF+b*F$.
Now $*$ is $SO(1,3)-$invariant, hence  
$SO(1,3)$ is a group of (complex) linear transformations of ${\cal G}_p$, 
preserving the quadratic form $\langle F,F \rangle=
-*\left(F\wedge (*F)+\IM  F\wedge F\right)$: this means
that this 'norm' is invariant by Lorentz transformations, 
hence it is of relevant physical
interest.
If we introduce the following coordinate functions
on ${\cal G}_p$:
$
z^1=F_{01}-iF_{23}
$,
$
z^2=F_{02}+iF_{13}
$ and
$
z^1=F_{03}-iF_{12}
$,
we have that
$
\displaystyle
\langle F,F \rangle=
(z^1)^2+(z^2)^2+(z^3)^2,
\label{complexeuclid}
$
 hence there naturally
arises the so called complex-Euclidean metric
on $\CI^3$: on one hand, by changing coordinates
we are brought to a generic symmetric bilinear
form on $\CI^3$; on the other one there arise 'poles'
if we attempt to extend 
the above construction e.g. to $(\PI^1)^3$.
Now the idea of generalizing to the curved framework is 
quite natural: the reader is referred to section
\ref{blabla}. 
Our main concern will be warped products of Riemann surfaces:
let $\,{\cal U}_i\simeq \DI$ or
$\,{\cal U}_i\simeq \CI$, with
coordinate function $u^i$
and
metric $b_1(u^1)\,du^1\odot du^1$
or $f_i(u^i)\,du^i\odot du^i$ if $i\geq 2$;
both $b_1$ and the $f_i$'s are nonzero meromorphic functions.
A {\sl warped product} of the $\{{\cal U}_i\}$'s
will be a meromorphic Riemannian manifold
(see definition \ref{hol_met})
$$
\left(
\prod_{i=1}^N {\cal U}_i,
b_1(u^1)\,du^i\odot du^i+
\sum_{i=2}^N a_i(u^1)f_i(u^i)
\,du^i\odot du^i
\right),
$$
where
the $a_k$'s ($k\geq 2$) are 
nonzero
meromorphic functions
(called {\sl warping functions})
defined on ${\cal U}_1$. 
This construction can be naturally generalized
to the case when the $\{{\cal U}_i\}$'s are more general Rieman surfaces.
We report that many of the known exact solutions of Einstein's 
field equations can be related, by means of 'complexifications', 
to such manifolds.

We introduce the concept of {\sl coercivity} of a warped product: informally speaking, it will amount to the fact that primitives
of 'square roots' of
some rational functions
of the coefficients involved in the metric
can be analitically continued until they 
take all complex values
but at most a finite number of ones.

Geodesics will show various types of 'singularities':
we record, among the other ones,
'logarithmic' singularities: they will be, more or less, points resembling $0$ in connection with 
$z\mapsto\log z$; rather more formally, 
a 'logarithmic singularity' $\ell$ will be a 
point 
in a two dimensional real
topological manifold, admitting 
a neighbourhooud ${\cal U}$ such that
${\cal U}\setminus\{\ell\}$
is a Riemann surface, but
 there is no complex structure
'at' $\ell$: this type of singularities arises from
the fact that geodesic equations admit first integrals whose solutions have poles with
nonzero residues.

We introduce the notion of {\sl completeness}:
a path will be essentially a holomorphic function
$F\colon S\rightarrow \M$, where $S$ is a Riemann surface over a region of $\PI^1$, admitting a projection mapping $\pi\colon
S\rightarrow\PI^1$: it will be {\sl complete}
provided that $\PI^1\setminus \pi(S)$ is a finite set: we are now able to attemtp to give a hazy idea
of 
our main result.
{\tt Theorem}: {\it a warped product of Riemann
surfaces is complete
(i.e. 'almost every' geodesic is complete)
 if and only if it is coercive.}

The last statement resumes the meaning of
theorems \ref{bla1}, \ref{bla2} and
\ref{bla3},
whilst definition of completeness is in \ref{completeness} and of coercivity
in \ref{coercive}.
We end this section with some references:
the problem of geodesic singularities arises from
semi-Riemannian geometry: see e.g.
\cite{beemerlich}; a different approach to holomorphic geometry could be found in \cite{manin}.
Finally, we owe 
\cite{lebrun} for the definition of 
a nondegenerate holomorphic metric, 
of a connexion (see p. 11 ff) and of a complex geodesic (see p.12 ff).
\section{Analytical continuation}
The idea of analytical continuation of a holomorphic mapping element 
$f:{\cal U}\rightarrow\M$ ( ${\cal U}$ is a region in the complex plane, $\M$, throughout this paper will be a complex manifold) is well known and amounts to a quintuple $Q_{\M}=(S,\pi,j,F,\M    )$, where
$S$ is a connected Riemann surface over a region of $\PI^1$,
$\pi\,\colon\, S\longrightarrow \CI$ is a nonconstant holomorphic mapping 
such that $U\subset \pi(S)$,
$j\,\colon\, U\longrightarrow S$ is a holomorphic  immersion such that $\pi\circ j=id\vert_{U}$
and 
$F\,\colon\, S\longrightarrow \M$ is a holomorphic mapping such that $F\circ j=f$.
Each finite branch point is 
kept into account by the fact of lying 'under' some critical point of $\pi$; it is a well known
(see e.g. \cite{cassa}, chap. 6) result that there exists a unique maximal analytical continuation,
called the {\BBB Riemann surface}, of $\left({\cal U},f   \right)$.
In the following we shall abbreviate 'holomorphic function element' by 'HFE' and 'holomorphic function germ' by 'HFG'.
For further purposes, we shall consider also 'poles' and 'logarithmic singularities': our definitions will axiomatize the behaviour of continuations of complex-valued holomorphic elements.
\begin{definition}
\rm
A {\BBB pole} of $Q_{\M}$ is a decreasing sequence of open sets $\{V_k\}_{k\geq K}\subset S$
such that there exist a positive integer $n$ and a point $ z_0\in\PI^1$,
such that $\bullet$ (P1)
for every $k\geq K$ $V_k$ is a connected component of 
${\pi^{-
1}(D(z_0,\frac{1}{k})\setminus \{z_0\})}$,
$\bullet $ (P2) for every $k\geq K$
${\pi\vert_{V_k}:V_k\longrightarrow
(D(z_0,\frac{1}{k})
\setminus \{z_0\})}$ is a n-sheeted 
covering and 
$\bullet$ (P3)
${\bigcap_{k\geq K}\overline{V_k}=
\emptyset}$
$\bullet$(P4)
there exist:
{an open set } $\Omega\subset\M$;
complex submanifolds $\N\subset\Omega$
{and }$\P\subset\Omega$ ($dim(\P)\geq 1$);
such that
$\Omega$ {and }$\N\times\P$  are
biholomorphic;
{for every } $k$, $ F(V_k\setminus\{p\}   )\subset\Omega$;
$pr_1\circ F:V_k\longrightarrow\N$ 
{has a removable singularity at } $p$ and
$\bigcap_{k\geq K}\overline{pr_2\circ F(V_k)}=
\emptyset$; a
{\BBB  logarithmic singularity} (in the following: {\BBB L-singularity}) $q$ of 
$Q_{\M}$  is a sequence of decreasing open sets $\{V_k\}_{K\geq K}$ of $S$ such that there
hold (P1), (P3) and
$\bullet$ (LS2)
for every $k\geq K$ and every (real) nonconstant closed path $\gamma:[0,1]\longrightarrow D(z_0,1/k)\setminus\{z_0\}$, with nonzero winding number around $z_0$, every lifted path $\beta:[0,1]\longrightarrow \pi^{-1}(D(z_0,1/k)\setminus\{z_0\})$ with respect to the topological covering $\pi$ is not a closed path, i.e. $\beta(0)\not=\beta(1)$; 
$q$ is 
$\bullet$ (RMLS)
a {\BBB removable L-singularity} for  if there exists $\eta\in\M$
such that $\bigcap_{k}\overline{ F(V_k)}=\{\eta\}$;
$\bullet$ (PLS)
a {\BBB polar L-singularity} for $\widehat F$ if there exist:
an open set $\Omega\subset\M$;
complex submanifolds $\N\subset\Omega$
{and }$\P\subset\Omega$ ($dim(\P)\geq 1$
)
such that
$\Omega$ and $\N\times\P$ are 
biholomorphic;
{for every } $k$, $F(V_k\setminus\{p\}   )\subset\Omega$;
$pr_1\circ F:V_k\longrightarrow\N$ 
{has a removable singularity at } $p$;
$\bigcap_{k=\geq K}\overline{pr_2\circ 
F(V_k)}=
\emptyset$;
\label{log_singol2}
\label{log_sing}
\label{sps}
\end{definition}

It is easily seen that $\{${\TTT L-singularities}$\}$ $\bigcap$ $\{${\TTT poles} $\}$ 
$=\emptyset$ and that $\widetilde S\colon=
S\cup$ $\{${\TTT poles of $Q_{\M}$} $\}$ has a canonical structure of a Riemann surface and $\pi$ admits a holomorphic extension $\widetilde\pi $ to $\widetilde S$, hence
an {\BBB extended analytical continuation} of $(U,f)$ is a quintuple $\widetilde Q_{\M}=(\widetilde S,\widetilde \pi,\widetilde j, F,\M    )$,
where
$\widetilde S$ and $\widetilde\pi$ are as above and
$\widetilde j=id_{S\rightarrow\tilde S}
\circ j$; of course there exists a unique maximal extended continuation of $\left({\cal U},f   \right)$,  build up as above, starting from its unique maximal continuation.

Consider now the set $B$ of the L-singularities of $Q_{\M}$: set $S^{\sharp}=S\bigcup B$ as a set and 
introduce a topology on $S^{\sharp}$:
open sets are the open sets in $S$ and a fundamental neighbourhood system of
the L-singularity $q={\{V_k\}_{k\geq K}}\in B$ is yielded by the sets $V_k^{\sharp}=V_k\bigcup \{q\}$.
\begin{lemma}\rm
$S^{\sharp}$ admits no complex structure at
 $q={\{V_k\}_{k\geq K}}$.
\end{lemma}
{\bf Proof:} were there one, we could find charts $({\cal W},\phi   )$ around $q$ and $({\cal V},\psi   )$ around $z_0$ such that 
$
\psi\circ\pi\circ\phi^{-1}(\zeta   )
=\zeta^N
$
for some integer $N>0$.
This fact would imply $\pi\vert_{{\cal W}\setminus\{q\}}$ to be a n-sheeted covering of ${\cal V}\setminus\{z_0\}$;
it is easily seen tha this fact would contradict
(LS2) in definition \ref{log_sing}.
\QUAN
\begin{lemma}\rm
{\bf (A)}: $\pi$ admits a unique continuous extension $\pi^{\sharp}$ to $S^{\sharp}$;
{\bf (B)}:  for every removable logarithmic 
singularity $r$ of $Q_{\M}$, $F$ admits a unique continuous extension
 $F^{\sharp}$ to $r$.
\label{funest}
\label{proest}
\end{lemma}
{\bf Proof:} {\bf (A)}: let $b\in B$ and $\{V_k\}$ be the sequence spotting $b$: define
$
\pi^{\sharp}(q)=\pi(q)$ if $q\in V_k$
and 
$
\pi^{\sharp}(b)=z_0 
$,
where $z_0$ is the common centre of the discs onto which the $V_k's$ are projected.
Now $\pi^{\sharp}$ is continuous at all points in $V_k$; moreover, for every neighbourhood $G$ of $z_0$, 
$\pi^{\sharp\ -1}(G)\supset \pi^{\sharp\ -1}(z_0) \bigcup\pi^{-1}(G\setminus\{z_0\})$,
hence, if we set  $H=\{b\}\bigcup\pi^{-1}(G\setminus\{z_0\})$, we have that $H$ is a neighbourhood of $b$ in $S^{\sharp}$ such that $\pi^{\sharp}(H)\subset G$, proving continuity at $b$. Arguing by density, we conclude that this extension is unique; the proof of {\bf (B)} is analogous.
\QUAN
\begin{definition}
\rm
A quintuple $ Q_{\M}^{\natural}=(S^{\natural},\pi^{\natural},j^{\natural},F^{\natural},\M)$, is 
an {\BBB analytical continuation with L-singularities} of the function element $(U,f)$ if there exists an analytical continuation
$ Q_{\M}$ of 
$(U,f)$ such that
$ S^{\natural}\setminus S$ consists of L-singularities of $F$, $\pi^{\natural}$ is the unique continuous extension of $\pi$ to
 $S^{\natural}$,
$j^{\natural}=id_{S\longrightarrow S^{\natural}}\circ j$ and 
$F$ admits a unique continuous 
extension $ F^{\natural}$ to $S^{\natural}\setminus \{ polar\ logarithmic\ 
singularities\ of \ F\}$.
$Q_{\M}^{\natural}$ is: {\BBB maximal} provided that so is $Q_{\M}$ and $Q_{\M}^{\natural}\setminus Q_{\M}$ contains all L-singularities of $Q_{\M}$; {\BBB extended} provided that so is $Q_{\M}$.
\label{loganalcont}
\label{ext_rsb}
\end{definition}
\begin{lemma}
\label{inverse}
\rm
{\bf 1)}: let $\f$ and $\g$ be two holomorphic germs each one inverse of the other;
 let $(   R,\pi,j,F,\CI)$ and $(S,\rho,\ell,G,\CI  )$
be their respective Riemann surfaces: then $F(R)=\rho(S)$;
{\bf 2)}: \label{quasiinverse}
\rm
let 
$\f$, $\g$, $\h$
be three HFG's such that
$\f\circ \g=\h$.
Let $(   R,\pi,j,F,\CI)$ 
be the Riemann surface of $\f$, $(S,\rho,\ell,G,\CI   )$
the one of $\g$ and 
$(T,\sigma,m,H,\CI  )$
the Riemann surface with L-singularities of $\h$: then $F(R)\setminus(\PI^1\setminus(
\sigma(T)
   )   )\subset\rho(S)$.

\end{lemma}
{\bf Proof:} we shall prove only {\bf 1)}; 
{\bf 2)} is analogous.
a) $F(R)\subset\rho(S)$: let $\xi\in R$ and $F(\xi)=\eta$; there exist:
an open neighbourhood ${\cal U}_1$ of $\xi$;
open subsets ${\cal U}_2\subset\pi({\cal U}_1   )$ and ${\cal V}_2\subset F({\cal U}_1   )$ and 
a biholomorphic function $g_2:{\cal V}_2\longrightarrow{\cal U}_2$, with inverse function $f_2:{\cal U}_2\longrightarrow{\cal V}_2$
{such that}:
$({\cal U}_2,f_2   )$ and $({\cal U},f  )$ are connectible and so are
$({\cal V}_2,g_2   )$ and $({\cal V},g  )$.
By construction there hence exist two holomorphic immersions 
$\widetilde j:{\cal U}_2\longrightarrow R\hbox{ and }
\widetilde\ell:{\cal V}_2\longrightarrow S$
such that $\pi\circ\widetilde j=\id$ and $\rho\circ\widetilde\ell=\id$.
Let ${\cal V}_1=F(U)_1$ and 
$
\Sigma=\{(x,y   )\in{\cal U}_1\times{\cal V}_2: F(x)=y\}
$;
moreover let $J:{\cal V}_2\longrightarrow\Sigma$
be defined by setting $ J(v)=(\widetilde j\circ  g_2(v),v )$.
Then $ (\Sigma,pr_2,J,\pi\circ pr_1   )$ is an analytical continuation of $({\cal V}_2,g_2   )$; indeed $ \pi\circ pr_1\circ J=\pi\circ\widetilde j\circ g_2=g_2$. But $({\cal V_2},g_2   )$ is connectible with $({\cal V},g     )$, hence
$ (\Sigma,pr_2,J,\pi\circ pr_1   )$ is an analytical continuation of $({\cal V},g )$.
There eventually exists a holomorphic function $h:\Sigma\longrightarrow S$ such that $\rho\circ h=pr_2$: hence 
$
\eta=pr_2(\xi,\eta   )=\rho\circ h
(\xi,\eta   )\in\rho(S   )
$.
b) $\rho(S)\subset F(R)$: let $s\in S$: there is a neighbourhood  $V$ of $s$ in $S$ such that $V\setminus\{s\}$ consists entirely of regular points both of $\rho$ and $G$, not excluding that $s$ itself be regular for $\rho$ or $G$ or both.
This fact means that for each $s^{\prime}\in V\setminus\{s\}$ there exists a HFE $({\rho(s^{\prime})},{\cal V}^{\prime},\widetilde g_{s^{\prime}}   )$ connectible with $({\cal V},g   )$ and, besides, a holomorphic immersion $\widetilde\ell:{\cal V}^{\prime}\longrightarrow V$.
By a) already proved, $G(s)\in\pi(R)$, hence there exist
$p\in R$ such that $\pi(p)=G(s)$ and
a neighbourhood $W$ of $p$ in $R$ such that $ \pi^{-1}(\widetilde g({\cal V}^{\prime}   )   )\bigcap W\not=\emptyset$.
Set
$ W^{\prime}=\pi^{-1}(\widetilde g({\cal V}^{\prime}   )   )\bigcap W$:
we may suppose, without loss of generality, that $\pi$ is invertible on $W^{\prime}$: hence there exists a (open) holomorphic immersion $\widetilde j:\widetilde g({\cal V}^{\prime}   )\longrightarrow W$.
Therefore, for each $\zeta\in\widetilde j(\widetilde g({\cal V}^{\prime}   )   )$, there exists $ \eta\in\widetilde\ell({\cal V}^{\prime}   )$ such that $ F(\zeta)=F(\widetilde j \circ \widetilde g\circ \rho(\eta)  )$.
Now, by definition of analytical continuation there holds $ F\circ\widetilde j\circ\widetilde g=\id$, hence we have $ F(\zeta)=\rho(\eta)$.
Consider now the holomorphic function 
$\Xi:W\times V\longrightarrow\CI$ 
defined by setting $ \Xi(w,v   )=F(w)-\rho(v)$: we have
$
(\Xi\vert_{\widetilde j(\widetilde g({\cal V}^{\prime}   )   )\times\widetilde\ell({\cal V}^{\prime})}
\equiv 0
$, but $ {\widetilde j(\widetilde g({\cal V}^{\prime}   )   )\times\widetilde\ell({\cal V}^{\prime})}$ is an open set in $W\times V$, hence $\Xi\equiv 0$ on $W\times V$, which in turn implies $ F(p)=\rho(s)$.
Therefore we have proved that for each $s\in S$ there exists $p\in R$ such that $ F(p)=\rho(s)$: this eventually implies that $\rho(S)\subset F(R)$.
\QUAN

\section{Complex-Riemannian metric 
structures}
\label{blabla}
\begin{definition}\rm
Let ${\cal E}$ be a closed hypersurface in $\M$:
{\sl an ${\cal E}$-meromorphic section} of 
${\cal T}_{r}^{s}\N$ is a 
holomorphic section $\Lambda$ of ${\cal T}_{r}^{s}(\M\setminus{\cal E})$ 
such that
for every $p\in{\cal E}$ and every 
chart
$\left({\cal U},(z^1...z^n)\right)$ around
$p$, there exists a neighbourhood 
$U$ of $p$ and $r\cdot s$ pairs of $\CI-$valued 
holomorphic functions $\phi_{i_1...i_r},\
\psi_{l_1...l_s}$, with $\psi_{l_1...l_s}\not=0$ on 
$U\setminus{\cal E}$, such that
$
\Lambda
\left( dz^{l_1}...dz^{l_s},
\frac
{\partial}
{\partial z^{i_1}}...
\frac{\partial}
{\partial z^{i_r}}
\right)=
\frac{\phi_{i_1...i_r}}
{\psi_{l_1...l_s}}
$.
{\TTT A 
complex
metric on} $\M$ is a 
symmetric section 
of ${\cal T}_{0}^{2}\M$.
It will be called {\sl
 holomorphic or ${\cal E}$-meromorphic 
}
provided that so is as a section;
$\Lambda$ is {\BBB nondegenerate} at $p$ if $rk(\Lambda(p))=dim(\M)$, {\BBB degenerate} otherwise;
if ${\cal D}$ is a closed hypersurface in 
$\M$ and $\Lambda$ is degenerate only 
on $\cal D$, we shall say that $\Lambda$ is 
${\cal D}$-degenerate.
We say that $p$ is a {\BBB metrically ordinary point} in $\M$ if $\Lambda$ is holomorphic 
and nondegenerate at $p$.
\label{reg_point}
\label{rank_met}\label{hol_met}
A {\BBB holomorphic 
(resp. nondegenerate holomorphic, resp. meromorphic)
Riemannian manifold} is a complex manifold endowed with a 
holomorphic (resp. nondegenerate holomorphic, resp. meromorphic)
metric.
\label{riemann}
\end{definition}
We now turn to introducing the holomorphic
Levi-Civita connexion induced by a
meromorphic, possibly degenerating metric.
First we need to introduce the holomorphic Levi Civita connexion induced by a holomorphic nondegenerate metric: this is done in a quite 
natural way.
Things are different if we allow metrics to
be meromorphic behaviour or to lower in their ranks.
These metric 'singularities' will be generally supposed to lie in closed hypersurfaces; Levi Civita  connexions may still be defined, but, as one could expect, they will turn out to be themselves 'meromorphic'.
Let now $(\N,\Lambda    )$ be a meromorphic Riemannian manifold admitting closed hypersurfaces ${\cal D}$ and ${\cal E}$ such that 
$\Lambda\vert_{\N\setminus{\cal E}}$ is holomorphic and
$\Lambda\vert_{(\N\setminus{\cal E})\setminus{\cal D}}$ is nondegenerate.
Since $\N\setminus{\cal E}$ is connected, we have that $(\N\setminus{\cal E})\setminus{\cal D},\Lambda\vert_{(\N\setminus{\cal E})\setminus{\cal D}}$ is a nondegenerate holomorphic Riemannian manifold admitting, as such, a canonical holomorphic Levi-Civita connexion $D$.
Now, if $p\in{\cal D}\bigcup{\cal E}$ and $V,W$ are holomorphic vector fields in a neighbourhood ${\cal V} $ of $p$ we 
can
define the vector field $D_V W$ on ${\cal V}\setminus({\cal D}\bigcup{\cal E}    )$, and this will be a ${\cal D}\bigcup{\cal E}-$meromorphic vector field.
The {\BBB Christoffel symbols } of a coordinate system $Z=(z^1\cdots z^m)$ on an open set ${\cal U}\subset\N$ are those complex valued functions, defined on
${\cal U}\setminus({\cal D}\bigcup{\cal E}    )$ by setting $\Gamma_{ij}^k =dz^k(  D_{\frac{\partial}{\partial z^i}}({\frac{\partial}{\partial z^j}}   )  )$.
Now the representative matrix $(g_{ij})$ of $\Lambda$ with respect to the coordinate system $Z$ 
is holomorphic in ${\cal U}$, with nonvanishing determinant function on 
${\cal U}\setminus({\cal D}\bigcup{\cal E}    )$; as such it admits a inverse 
matrix ${g^{ij}}$, whose coefficients hence result in being ${\cal D}\bigcup{\cal E} 
$-meromorphic functions.
It is easy to prove that
$D_{\frac{\partial}{\partial z^i}}(\sum_{j=1}^m W^j {\frac{\partial}{\partial z^j} }   )=
\sum_{k=1}^m(\frac{\partial W^k}{\partial z^i}+\sum_{j=1}^m \Gamma_{ij}^k W^j   )\frac{\partial}{\partial z^k}$ as meromorphic vector fields and
$2\Gamma_{ij}^k =\sum_{m=1}^N g^{km}(-g_{ij,m}+g_{im,j}+g_{jm,i})=2\Gamma_{ij}^k $ as meromorphic functions; then:
\begin{proposition}\rm
For every pair $V,W$ of holomorphic vector field on the open set ${\cal U}$ (belonging to a maximal atlas) in the meromorphic Riemannian manifold $(\N,\Lambda     )$, $D_V W$ is a well defined vector field, holomorphic on ${\cal U}\bigcap\{n\in\N:\Lambda\hbox{ is holomorphic and nondegenerate at $n$}\}$ and may be extended to a meromorphic vector field on ${\cal U}$.
\end{proposition}
{\bf Proof:} there exist holomorphic functions $\{V^i\}$, $\{W^j\}$ and a coordinate system $Z=(z^1.....z^N    )$
on ${\cal U}$ such that 
$
{V=\sum_{i=1}^N V^i\frac{\partial}
{\partial z^i}}$ and
$ {W=\sum_{j=1}^N W^i\frac{\partial}
{\partial z^j}}
$.
The fact that
$
D_V W=\sum_{i=1}^N V^i D_{\frac{\partial}{\partial z^i}}
(\sum_{j=1}^N W^j   \frac{\partial}{\partial z^i} )=
\sum_{k=1}^N (\sum_{i,j=1}^N V^i( \frac{\partial W^k}
{\partial z^i}+\Gamma^k_{ij}W^j   )    )\frac{\partial}{\partial z^k}
$
ends the proof.\QUAN
\begin{definition}
\rm
Given a ${\cal D}$-degenerate and ${\cal E}$-meromorphic Riemannian manifold $(\N,\Lambda    )$, with ${\cal
D} $ and ${\cal E}$ closed hypersurfaces in $\N$, the {\BBB Levi-Civita metric connexion} (or {\BBB meromorphic metric
connexion}) of $\N$ is the collection consisting of all metric connexions $\{D\left[{\cal U}_i\setminus({\cal D}\bigcup
{\cal E}  ) \right]\}_{i\in I}$ as ${\cal U}\}_i $ runs over any maximal atlas ${\cal B}=(\{{\cal U}\}_i    )_{i\in I}$ on $\N$. 
\end{definition}

\subsection{Meromorphic parallel translation
and geodesics}
We now slightly reformulate the notion of path to cope with the complex environment:
a {\BBB path} in $\M$ is a quintuple $Q_{\M}=(S\bigcup,\pi,j,F,\M)$,
where
$S$ is a connected Riemann surface, 
$\pi\in{\cal H}\left(S,\PI^1   \right)$, 
$F\in {\cal H}\left(S,\M   \right)$ and $j$
is a holomorphic immersion $j\colon U\longrightarrow S\setminus\Sigma$ such that $\pi\circ j=id\vert_{U}$, where $U$ is a region in the complex plane; a path is $z_0-$starting at $m$ provided that $z_0\in U$ and $F\circ j(z_0)=m$.

In the continuation,
we shall call $T\M$ (resp.$T^{*}\M$) $\M$'s holomorphic tangent (resp. cotangent) bundle 
and, more generally, ${\cal T}_{r}^{s}\M$ its holomorphic $r-$covariant and $s-$contravariant tensor bundle; as usual, $\Pi\colon{\cal T}_{r}^{s}\M\longrightarrow\M$ will denote
their natural projections.
We now define the {\BBB velocity field} of a 
path $Q_{\M}$ as a suitable meromorphic section over $F$ of the holomorphic tangent bundle $T\M$: to achieve this purpose, we need to lift the vector field ${d}/{dz}$ on $\CI$ with respect to $\pi$; of course, in general, contravariant tensor fields couldn't be lifted, but
we may get through this obstruction by keeping into account that $\CI$ and $S$ are one-dimensional and allowing the lifted vector field to be meromorphic.
We call $P$ the set of branch points of $\pi$.
\begin{lemma}\rm
there exists a unique $P$-meromorphic vector field $\widetilde{{d}/{dz}}$ on $S$ 
such that, for every $r\in S\setminus P$, $\pi_{*}\vert_r(\widetilde{{d}/{dz}}\vert_{r}    )=({d}/{dz})\vert_{\pi(r)}$.
\label{lift}
\end{lemma}
{\bf Proof:} consider $\omega=\pi^{*}dz$ and  $\Lambda=\pi^{*}(dz\odot dz    )$ on $S$: the latter establishes an isomorphism between the holomorphic cotangent  and  tangent bundles
of $S\setminus P$.
Call $V$ the holomorphic vector field corresponding to $\omega$ in the above isomorphism: we claim that $V=\widetilde{{d}/{dz}}$ on $S\setminus P$.
To show this fact, we explicitely compute the components of $V$ with respect to
a maximal atlas ${\cal B}=\left\{(U_{\nu},\zeta_{\nu}    )\right\}$
for $S\setminus P$: let
$
\omega_{(\nu)\ 1}=
\omega(\partial/\partial \zeta_{(\nu)}    )$,
$
g_{(\nu)\ 11}=\Lambda(\partial/\partial \zeta_{(\nu)},\partial/\partial \zeta_{(\nu)}    )
$;
then, set ${V_{(\nu)}^1={\omega_{(\nu)\ 1}}/{g_{(\nu)\ 11}}}$ the collection $\left\{( {\cal U_{\nu}}, V_{(\nu)}^1  )\right\}$ of open sets and holomorphic functions is such that,
on overlapping local charts $(U_a,\zeta_a    )$ and $(U_b,\zeta_b   )$,
we have
$$
V_{(a)}^1=
\frac
{\omega_{(a)\ 1}}{g_{(a)\ 11}}
=\frac
{\omega_{(b)\ 1} (d\zeta_{(b)}/d\zeta_{(a)})}
{g_{(b)\ 11}({ d\zeta_{(b)}}/{d\zeta_{(a)}})^2}
=
{{V_{(b)}^1}
\frac
{d\zeta_{(a)}}{d\zeta_{(b)}}} ,
$$
that is to say that collection defines
a holomorphic
vector field. Now for every $r\in S\setminus P$, 
$$
dz\vert_{\pi(r)}
( 
\pi_{*}\vert_r
\widetilde{{d}/{dz}}\vert_{r}
)=\pi^*dz\vert_r
(
\widetilde{{d}/{dz}}\vert_{r}
)
=\frac{\pi^*dz\vert_r (\partial/\partial\zeta\vert_r) }
{dz\vert_{\pi(r)}(\pi_*\partial/\partial\zeta\vert_r) 
}
=1
,$$
hence
$
\pi_{*}\vert_r(\widetilde{{d}/{dz}}\vert_{r}    )=({d}/{dz})\vert_{\pi(r)}
$,
proving the asserted.
 
Let's prove that $\widetilde{d/dz}$ may be extended to a meromorphic vector field on $S$:
if $p\in P$ then
we can find local charts $(U,\psi    )$ around $p$, $(V,\phi)$ around $\pi(p)$,
and an integer $N>0$ such that $\phi\circ\pi\circ\psi^{-1}(u)=u^N$.
Now we have
$$
(\psi^{-1\ *}\pi^{*}\phi^{*}(dw)\frac d {du})(u)=
dw(\phi_{*}\pi_{*}\psi_{*}^{-1}\frac d {du})\vert_u)
=dw((\phi\pi\psi^{-1})^{\prime}\frac d {dw}))=Nu^{N-1};
$$
 but $\phi$ and $\psi$ are charts, hence $\pi^{*}dz $ itself is vanishing of order $N-1$ at $p$;
as already proved, $\pi_{*}\vert_r(\widetilde{{d}/{dz}}\vert_{r}    )=({d}/{dz})\vert_{\pi(r)}$ on $U\setminus \{p\}$ and, consequently,
$(\pi^{*}dz)(\widetilde{d/dz})=dz(\pi_{*}\widetilde)=dz(d/dz)=1$ on $U\setminus \{p\}$, 
hence on $U$.
Now, in local coordinates, 
$
(\pi^{*}dz)=\alpha d\phi$ and
$\widetilde{d/dz}=y\,{\partial}/{\partial\phi}
$,
where $\alpha$ is a holomorphic function on $U$, vanishing of order $N-1$ at $p$ and $y$ is a holomorphic function on $U\setminus \{p\}$.
By the argument above, $y\alpha=1$, hence $y$ has a pole of order $N-1$ at $p$:
a similar argument holds for each isolated point in $P$, 
proving the meromorphic behaviour of $\widetilde{d/dz}$.
\QUAN
\begin{lemma}\rm
The mapping 
$
r\longmapsto 
(F,F_{*}\vert_r
(\widetilde{\frac{d}{dz}}\vert_r
)
)$
may be extended to a P-meromorphic section  of 
$T\M$ over $F$.
\label{velo_def}
\end{lemma}
{\bf Proof:}
let $p\in P$ and $U$ be a neighbourhood of $p$ such that
there exist
a local chart 
$\zeta : U\longrightarrow\CI_w$ and holomorphic functions 
$f,g$ on $\zeta(U)$  such that 
$\widetilde{\frac{d}{dz}}\vert_{\zeta_{-1}(U)}=
\zeta^{-1}_{*}(\frac{f}{g}(w)\frac{d}{dw}\vert_w)$;
for every local chart $\Psi=(u^{1...m},du^{1...m})$ on $T\M$,
\begin{eqnarray*}
\Psi\circ V \circ \zeta^{-1}(w)&=&
\Psi\circ (
F\circ \zeta^{-1}(w),F_{*}\vert_{\zeta^{-1}(w)}\zeta^{-1}_{*}(\frac{f}{g}(w)
\frac{d}{dw}\vert_w))\\
&=&\Psi\circ ( F\circ \zeta^{-1}(w),\frac{f}{g}(w)\frac{d}{dw}(F\circ\zeta^{-1})(w))\\
&=&
\left(u^{1...m}\circ F\zeta^{-1}(w), 
\frac{f}{g}(w)\frac{d}{dw}(u^{1...m}\circ F\circ\zeta^{-1})(w)\right)
\end{eqnarray*}
\QUAN
\begin{definition}\rm
The {\BBB velocity field } of a path $Q_{\M}=(S,\pi,j,F,\M)$ is the meromorphic mapping
$V(Q_{\M})\colon S\setminus P \longrightarrow T\M$ defined by
$ r\longmapsto (F,F_{*}\vert_r
(\widetilde{\frac{d}{dz}}\vert_r)
)$
\label{vf}
\end{definition}

We turn now to study vector fields on paths: an obvious example is the velocity field, defined in definition \ref{vf}: just as in semi-Riemannian geometry, there is a natural way of defining the rate of change $X^{\prime}$ of a meromorphic vector field $X$ on a path.
We study at first paths with values in a nondegenerate holomorphic Riemannian manifold $\M$:
let
$Q_{\M}=(S,\pi,j,\gamma,\M)$ be a path in $\M$;
$P$ be the set of branch points of $\pi$;
$r\in S\setminus P$ be 
such that $\widetilde{d/dz}$ is holomorphic at $r$, ${\cal V}\subset S\setminus P$
 be a neighbourhood of $r$ such that $\gamma({\cal V})$ is contained in a local chart in $\M$;
${\cal H(V    )}$ be the ring of holomorphic functions on ${\cal V}$ and 
${\cal X}_{\gamma}({\cal V}    )$ the Lie algebra of holomorphic vector fields 
over $\gamma$ on ${\cal V}$: it is well known that
there exists a unique mapping $\nabla_{\gamma^{\prime}}\colon
{\cal X}_{\gamma}({\cal V}    ) \longrightarrow {\cal X}_{\gamma}({\cal V}    )$,
called {\BBB induced covariant derivative} on $Q_{\M}$ such that
$\nabla_{\gamma^{\prime}}(aZ_1+bZ_2    )=a\nabla_{\gamma^{\prime}}
Z_1+b\nabla_{\gamma^{\prime}} Z_2$,
$\nabla_{\gamma^{\prime}} (hZ)=(\widetilde{\frac{d}{dz}}h    )Z
+h\nabla_{\gamma^{\prime}} Z,\quad h\in{\cal H(V    )}$
and $
\nabla_{\gamma^{\prime}} (V\circ\gamma )(r)=
D_{\gamma_{*}\vert_r(\widetilde{\frac{d}{dz}}\vert_r)}$,
where $V$
is a holomorphic vector field in a neighbourhood of $\gamma(r)$.
Moreover,
$\widetilde{\frac{d}{dz}}\left\langle X,Y\right\rangle=\left\langle \nabla_{\gamma^{\prime}}
 X,Y   \right\rangle+\left\langle X,\nabla_{\gamma^{\prime}} Y \right\rangle
\quad X,Y\in {\cal X}_{\gamma}({\cal V})$.
Now let ${\cal R}=\{{\cal V}_k\}_{k\in K}$ be a maximal atlas for $S\setminus P$; we may assume that, for every $k$,
maybe shrinking ${\cal V}_k   $, $\gamma( {\cal V}_k   )$ is contained in some local chart ${\cal U}_i$ in the
already introduced atlas ${\cal A}$ for $\M$.

Now, if ${\cal V}_1  $ and ${\cal V}_2  $ are overlapping open sets in ${\cal R}$, ${\cal V}_1 \bigcap{\cal V}_2\in {\cal R} $ too, and 
$
\nabla_{\gamma^{\prime}}
\left[ {\cal V}_1   \right]
\vert_{{\cal V}_1 \bigcap{\cal V}_2}
=\nabla_{\gamma^{\prime}}
\left[ {\cal V}_2   \right]
\vert_{{\cal V}_1 \bigcap{\cal V}_2}
$.
Now let's complete ${\cal R}$ to an atlas ${\cal S}$ for $S$: keeping into account that the local coordinate expression of the induced covariant derivative is
$
\nabla_{\gamma^{\prime}}Z=\sum_{k=1}^m( \widetilde{\frac{d}{dz}}Z^k+\sum_{i,j=1}^m\Gamma_{ij}^k
\widetilde{\frac{d}{dz}}(u^i\circ\gamma)Z^j   )\frac{\partial}{\partial u^k}
$,
hence  pairs of holomorphic vector fields on $\gamma$ are transormed into $P$-meromorphic vector fields on $\gamma$.
\begin{definition}\rm
The $P$-meromorphic {\BBB induced covariant derivative}, or the $P$-meromorphic parallel translation on a path $Q_{\M}=(S,\pi,j,\gamma,\M    )$ with set of branch points $P$ and taking values in a nondegenerate Riemannian manifold $\M$ is the collection consisting of the induced covariant derivatives $\nabla_{\gamma^{\prime}}\left[ {\cal V}_k \setminus P  \right] $ as ${\cal V}_k $ runs over a maximal atlas ${\cal S}=(\{{\cal V}_k\}    )_{k\in K}$ on $S$. 
\label{indcovdev}
\end{definition}
Let's turn now to dealing with meromorphic parallel translations induced on a path $Q_{\N}=(T,\varrho,j,\delta,\N    )$, in a meromorphic Riemannian manifold
$(\N,\Lambda)$ admitting closed hypersurfaces ${\cal D}$ and ${\cal E}$ such that 
$\Lambda\vert_{\N\setminus{\cal E}}$ is holomorphic and
$\Lambda\vert_{(\N\setminus{\cal E})\setminus{\cal D}}$ is nondegenerate.
We set ${\cal F}={\cal D}\bigcup{\cal E}$ and restrict our attention to paths $z_0$-starting at metrically ordinary points, supposing, without loss in generality, that $z_0=0$.
\begin{lemma}\rm
Set $\M=\N\setminus{\cal F}$, $S=\delta^{-1}(\M)$: then $T\setminus S$ is discrete, hence $S$ is a connected Riemann surface.
\end{lemma}
{\bf Proof:} suppose that there exists a subset ${\cal V}\subset T\setminus S$ admitting an accumulation point $t\in{\cal V}$ and consider a countable atlas for ${\cal B}=\{U_n\}_{n\in\ENNE}$ for $\N$ such that, for every $n$, there exists $\Psi_n\in{\cal O}( \{U_n\}   )$ such that 
$
U_n\bigcap{\cal F}=\{X\in U_n : \Psi_n=0\}
$.
Set $\delta^{-1}(U_n)=T_n\subset T$ and suppose, without loss of generality, that $\delta(t)\in U_0$.
Now $\Psi_0\circ\delta\vert_{{\cal V}\cap T_0}=0$ and $t\in{\cal V}\cap T_0$ is an accumulation point of ${\cal V}\cap T_0$ , hence $\Psi_0\circ\delta\vert_{ T_0}=0$ and $\delta(T_0)\subset {\cal F}$.
Suppose now that $T_N\not=\emptyset$ for some $N$: we claim that this implies $\delta(T_N)\subset{\cal F}$: to prove the asserted, pick two points $\tau_0\in T_0$ and $\tau_n\in T_n$ and two neighbourhoods $T^{\prime}_0$, $T^{\prime}_N$ of $\tau_0 $ and $\tau_n$ in $T_0$ and $T_n$ respectively, such that $\varrho\vert_{T^{\prime}_0}$ and $\varrho\vert_{T^{\prime}_N}$ are biholomorphic functions.
Now the function elements $(\varrho(T^{\prime}_0),\delta\circ( \varrho\vert_{T^{\prime}_0}   )^{-1}    )$ and $(\varrho(T^{\prime}_N),\delta\circ( \varrho\vert_{T^{\prime}_N}   )^{-1}    )$ are connectible, hence there exists a finite chain $\{W_{\nu}\}_{\nu=0...L}$ such that $W_0=\varrho(T^{\prime}_0)$, $W_L=\varrho(T^{\prime}_N)$, $W_{\nu}\bigcap W_{\nu+1}\not=0$ for every $\nu$.
Without loss of generality, we may suppose that each $W_{\nu}$ admits a holomorphic, hence open, immersion $j_{\nu}\longrightarrow T$, hence, setting
$
S_0=T_0$,
 $S_{\lambda}=j_{\lambda}(W_{\lambda}) \hbox{\rm\  for }
\lambda=1...L$,
 $S_{L+1}=T_N$
yields a finite chain of open subsets $\{S_{\lambda}\}_{\lambda=0...M}$ of $T$ connecting $T_0$ and $T_N$.
Let's prove, by induction, that, for every $\lambda$, $\delta(S_{\lambda})\subset{\cal F}$.
$\bullet$ {\BBB At first} recall that
$\delta(S_0)\subset U_0\bigcap{\cal F}$ as already proved; suppose now that $\delta(S_{k-1})\subset {\cal F}$. 
We have $S_{k-1}\bigcap S_k\not=\emptyset$, hence $\delta(S_{k-1})\bigcap\delta(S_{k})\not=\emptyset$.
For every $m$ set $\Sigma_{km}=\delta(S_{k-1})\bigcap\delta(S_{k})\bigcap U_m$: if $\Sigma_{km}\not=\emptyset$, then $\Psi_m\circ\delta\vert_{\delta^{-1}(\Sigma_{km})\bigcap S_{k-1}\bigcap S_k }\equiv 0$; but $\delta^{-1}(\Sigma_{km})\bigcap S_{k-1}\bigcap S_k$ is open in $\delta^{-1}(\delta(S_k)\bigcap U_m    )\bigcap S_k$, thus
$\Psi_m\circ\delta\vert_{\delta^{-1}(\delta(S_k)\bigcap U_m    )\bigcap S_k}\equiv 0$, that is to say $\delta(S_k)\bigcap U_m\subset{\cal F}$.
$\bullet$ {\BBB On the other hand},
if $\Sigma_{km}=\emptyset$, but $\delta(S_k)\bigcup U_m\not=\emptyset$ we claim that $\delta(S_k)\bigcap U_m\subset{\cal F}$ as well: proving this requires a further induction: pick a $U_M$ such that $\Sigma_{kM}\not=\emptyset$ and a finite chain of open sets
${\cal B}^{\prime}=\{U^{\prime}_{\mu}\}_{\mu=0...J}\subset {\cal B}$ (with $U_{\mu}^{\prime}\bigcap\delta(S_k)\not=\emptyset$ for each $\mu$) connecting $U_M$ and $U_m$.
Since $\Sigma_{kM}\not=\emptyset$, 
$\delta(S_k)\bigcap U_0^{\prime}=\delta(S_k)\bigcap U_M\subset {\cal F}$; suppose by induction that 
$\delta(S_k)\bigcap U_{l-1}^{\prime}\subset {\cal F} $:
then
$\Psi_l\circ\delta\vert_{\delta^{-1}(\delta(S_k)\cap U_{l-1}^{\prime}\cap U_{l}^{\prime}   )\cap S_k}\equiv 0$,
hence
$
\Psi_l\circ\delta\vert_{\delta^{-1}(\delta(S_k)\cap U_{l}^{\prime}   )\cap S_k}\equiv 0$,
i.e. $\delta(S_k)\bigcap U_{l}^{\prime}  \subset {\cal F}$: this fact ends the induction and eventually implies 
$\delta(S_k)\bigcap U_{m}= \delta(S_k)\bigcap U_{J}^{\prime}\subset {\cal F}$.
Summing up, $\delta(S_k)=\bigcup_m ( \delta(S_k)\bigcap U_m   )\subset {\cal F}$, for each $k$;  hence $\delta(T_N)=\delta(S_M)\subset{\cal F}$ and eventually $\delta(T)=\delta(\bigcup_{N\in\ENNE} T_N    )\subset {\cal F}$, hence $\delta$ couldn't start at a point in $\N\setminus{\cal F}$.
\QUAN
In the following considerations, there will still hold all notations introduced in 
preceding lemma: given a path $Q_{\N}=(T,\varrho,j,\delta,\N    )$, set $\pi=\varrho\vert_S$,
$\gamma=\delta\vert_S$ and note that, since $Q_{\N}$ is starting from a metrically ordinary point $m$, $j$ may be
supposed to take values in fact in $S$; since the preceding lemma shows that $S$ is a connected Riemann surface,
$Q_{\M}=(S,\pi,j,\delta\vert_S\M    )$ is in fact a path in $\M$, which we call the {\BBB depolarization} of $Q_{\N}$. But $\M$ is a nondegenerate holomorphic Riemannian manifold, hence if $P$ is the set of branch points of $\pi$, there is a $P$-meromorphic induced parallel translation on $Q_{\M}$, built up as in definition \ref{indcovdev}.
Finally, we introduce a maximal atlas ${\cal T}$ for $T$ and yield the following:
\begin{definition}\rm
Let $(\N,\Lambda    )$ be a ${\cal E}$- meromorphic and ${\cal D}$-degenerate Riemannian manifold,
$\M=\N\setminus({\cal D}\bigcup{\cal E}    )$,
$Q_{\N}=(T,\varrho, j ,\delta,\N    )$ 
a path: the {\BBB $(P\bigcup\delta^{-1}( {\cal D}\bigcup{\cal E}   )    )$-meromorphic induced covariant derivative} on $Q_{\N}$
 is the collection consisting of all induced covariant derivatives $\nabla_{\gamma^{\prime}}\left[ {\cal V}_k \bigcap S
 \right] $ as ${\cal V}_k$ runs over a maximal atlas ${\cal T}=(\{{\cal V}_k\}    )_{k\in K}$ for  $T$ and
 $Q_{\M}=(S,\pi,j,\delta\vert_S\M    )$ is the depolarization of $Q_{\N}$.  
\label{indcovdev2}
A meromorphic (in particular, holomorphic) vector field $Z$ on a path is {\BBB parallel} provided that $\nabla Z=0$ (as a {\BBB meromorphic} vector field).
\label{parall}
A {\BBB geodesic} in a meromorphic (in particular, holomorphic) Riemannian manifold is a path whose (meormorphic) velocity field is parallel.\label{geodesic}
\end{definition}
The local equations ${ \BETA^{\bullet\bullet}{}^{k}+\sum_{i,j=1}^N \Gamma_{ij}^k(\beta)
\BETA^{\bullet}{}^{i}\BETA^{\bullet}{}^{j}=0,
k=1.....N
}
$
of elements of geodesics $(U,\beta   )$
are a system of $N$  second-order o.d.e.'s in the complex domain, 
with meromorphic coefficients, in turn equivalent to an autonomous system of $2N$ 
first-order equations; as a consequence of
general theory (see e.g. \cite{hille}, th. 2.2.2) 
for every metrically ordinary point $p\in\M$, every holomorphic tangent vector 
$V_p\in T_p\M$ and every $z_0\in\CI$, there exists a unique germ 
\hbox{\boldmath{}$\beta_{z_0}$\unboldmath}
of geodesic such that 
\hbox{\boldmath{}$\beta_{z_0}$\unboldmath $(z_0)=p$} and
\hbox{\boldmath{}$\beta_{z_0\ *}$\unboldmath$(d/dz)\vert_{z_0}=V_p$;}
moreover any continuation of \hbox{\boldmath{}$\beta_{z_0}$\unboldmath}
is a geodesic. 
\begin{definition}
\rm
\label{completeness}
A meromorphic Riemannian manifold is {\TTT  complete} provided that the Riemann surface, with L-singularities, of each geodesic
starting at a metrically ordinary point is complete
\end{definition}
\section{Completeness theorems}
In this section we shall be concerned with warped products of Riemann surfaces, 
each one endowed with some meromorphic metric: in this framework we shall prove a 
geodesic completeness criterion.
Consider at first, like in foreword, a warped product of unit discs
or complex planes
$
\left(
\prod_{i=1}^N {\cal U}_i,
b_1(u^1)\,du^i\odot du^i+
\sum_{i=2}^N a_i(u^1)f_i(u^i)
\,du^i\odot du^i
\right)
$: in the following we shall denote it by
$$
{\cal U}={\cal U}_1\times_{a_2(u^1)}
{\cal U}_2\times_{a_3(u^1)}
{\cal U}_3\times
........
\times_{a_N(u^1)}{\cal U}_N
,
$$
and call it 
a {\BBB direct manifold}.
We recall that 
$b$, the $a_k$'s and the $f_k$'s are 
nonzero meromorphic functions, with
 $b$ and the $a_k$'s
defined on ${\cal U}_1$.
Each element of geodesic of $({\cal U},\Lambda   )$ satisfies the following
system of $N$ o.d.e.'s in the complex domain: 
\begin{equation}
\label{equazioninormali}
\cases{
\U^{\bullet\bullet}{}^1(z)+
\frac{b_1^{\prime}(u^1(z))}{2b_1(u^1(z))}
(\U^{\bullet}{}^1(z)   )^2
-\sum_{l=2}^N \frac{a_l^{\prime}(u^1(z))f_l(u^l(z))}
{2b_1(u^1(z))}(\U^{\bullet}{}^l (z)  )^2=0\cr
\ \cr
\U^{\bullet\bullet}{}^k(z)+
\frac{f_k^{\prime}(u^k(z))}{2f_k(u^k(z))}(\U^{\bullet}{}^k(z))^2
+\frac{a_k^{\prime}(u^1(z))}{a_k(u^1(z))}
(\U^{\bullet}{}^k (z))(\U^{\bullet}{}^1 (z)  )
=0,\ k=2...N,
}
\end{equation}
provided that it starts at a metrically ordinary point.
\begin{lemma}\rm
\label{integraleprimo1}
The equations (\ref{equazioninormali})
admit the following first integrals:
\begin{equation}
\cases
{
(A)\ if\ u^1\not=const\ \ 
\cases{
(
\U^{\bullet}{}^1(z) 
)^2
( b_1(u^1(z)   )   )
=A_1-\sum_{l=2}^N\frac{A_l}
{a_l(u^1(z)   )}
\quad\spadesuit\cr
(\U^{\bullet}{}^k(z)   )^2
f_k(u^k(z)   )
\left[a_k(u^1(z)   )\right]^2
=A_k\quad k=2...N\ \clubsuit.
}
\cr
(B)\ if\ u^1=const\ \
\cases{
{u}^1(z)
=A_1\quad\diamondsuit\cr
(\U^{\bullet}{}^k(z)   )^2
f_k(u^k(z)   )
=A_k\quad k=2...N\quad  \heartsuit,
}
}
\end{equation}
where the $A_k$'s are suitable complex constants.
\label{integraleprimo2}
\end{lemma}
{\bf Proof:} we prove only $(A)$; 
$(B)$ is analogous.
Divide
the k-th equation in
(\ref{equazioninormali}) by $u^k$ and
integrate once: then
$$
(\U^{\bullet}{}^k(z)){}^2
f_k(u^k(z)   )
\left[a_k(u^1(z)   )\right]^2
=(\U^{\bullet}{}^k(z_0)   )^2
f_k(u^k(z_0)   )
\left[a_k(u^1(z_0)   )\right]^2
\colon =
A_k
.$$
As to $\,\spadesuit\,$,
by the first equation of (\ref{equazioninormali}) there holds 
$
2b_1(u^1(z)   )\U^{\bullet}{}^1(z)
\U^{\bullet\bullet}{}^1(z)+{b_1^{\prime}(u^1(z))}
(\U^{\bullet}{}^1(z)   )^3
-\sum_{l=2}^N {a_l^{\prime}(u^1(z))f_l(u^l(z))}
(\U^{\bullet}{}^l (z)  )^2\U^{\bullet}{}^1(z)=0
$;
by $\clubsuit$, already proved,
$
(\U^{\bullet}{}^{l} (z)  )^2 
{f_l(u^l(z))\left[a_l(u^1(z))\right
]^2}=
{A_l}
$,
hence
$$
b_1(u^1(z)   )\U^{\bullet}{}^1(z)
\U^{\bullet\bullet}{}^1(z)+{b_1^{\prime}(u^1(z))}
(\U^{\bullet}{}^1(z)   )^3
-\sum_{l=2}^N 
A_l
\frac
{a_l^{\prime}(u^1(z))}
{\left[a_l(u^1(z))\right]^2}
\U^{\bullet}{}^1(z)=0;
$$
integrating once,
dividing by $b_1(u^1(z)   )$ 
and 
setting
$ A_1=K/b_1(u^1(z_0)   )$
ends the proof.
\QUAN
\begin{definition}
\rm
\label{coercive}
A direct manifold
$
{\cal U}
$
with metric
$
\Lambda(u^1.....u^N)=
b_1(u^1)\,du^i\odot du^i+
\sum_{i=2}^N a_i(u^i)f_i(u^i)
\,du^i\odot du^i
$,
where $b_1$, the $a_k$'s and the $f_k$'s are 
nonzero meromorphic functions
is {\BBB coercive} provided that, for every metrically ordinary
 point  $\displaystyle X_0=\left(x_0^1...x_0^N \right)$ and

$\bullet$ for every n-tuple 
$ (A_1...A_N   )\in\CI^N$ such that 
$
b_1(x_0^1)\not=0$, $A_1-\sum_{l=2}^N\frac
{
A_l}
{
a_l(x_0^1)}\not=0
$
and, for each one of the two HFG's 
${\alph_1}$
and
${\alph_2}$,
such that 
$
({\alph_i})^2=
\left[
\frac{1}{b_1}
(
A_1-\sum_{l=2}^N\frac
{A_l}
{a_l}
)
\right]
_{0}
\quad i=1,2
$,
the Riemann surface $(S_1,\pi_1,j_1,\Phi_1,{\cal U}   )$
of both the HFG's
$\left[
\int_{x_0}^{u^1}
{
\frac{d\,\eta}
{\alph_i(\eta)}
}
\right]_{x_0^1}
\ i=1,2;
$
is such that $\PI^1\setminus\Phi_1(S_1)$ is a finite set;

$\bullet$ 
for each $k$, $2\leq k\leq N$
and for each one of the two HFG's 
${\phi_{k1}}$
and
${\phi_{k2}}$
such that 
$
({\phi_{ki}}   )^2=
\left[
f_k
\right]_{x_0^1},
\ i=1,2
$,
the Riemann surface $(S_k,\pi_k,j_k,\Phi_k,{\cal U}   )$
of both the HFG's
$
\left[\int_{x_0^1}^{u^k}
\phi_{ki}(\eta)
\,d\eta\right]_{x_0^1}\ \ i=1,2
$
is such that $\PI^1\setminus\Phi_k(S_k)$ is a finite set.
\end{definition}
{\BBB Definition \ref{coercive} may be checked
for just one metrically ordinary point $X_0$: this is proved in lemma \ref{isom}; moreover, we may assume,without loss of generality 
$X_0=0$}: 
were not, we could carry it into $0$ by
applying an automorphism of 
${\cal U}$, that is to say a direct product of automorphisms of the unit ball or of the complex plane, according to the nature of each ${\cal U}_i$. 
Then a simple pullback procedure would yield back the initial situation: {\BBB in the following we shall understand this choice}.

In the following lemma we shall use the 'square root' symbol in the meaning of definition \ref{coercive}: in other words, given a HFG, which is not vanishing at some point, it should denote any one of the two HFG's yielding it back when squared.
\begin{lemma}
\label{isom}
\rm
For every metrically ordinary point $(\xi^1...\xi^N   )$ of ${\cal U}$ and every n-tuple $(A_1...A_N   )\in\CI^N$
such that
$
b_1(x_0^1)\not=0,\ 
A_1-\sum_{l=2}^N\frac
{
A_l}
{
a_l(x_0^1)}\not=0,\ 
\ b_1(\xi^1)\not=0,\ 
A_1-\sum_{l=2}^N\frac
{
A_l}
{
a_l(\xi^1)}\not=0
$,
set $\Psi(\eta)={ A_1-\sum_{l=2}^N\frac
{A_l}
{a_l(\eta)}}$,
the Riemann surfaces of the HFG's
$
\left[
\int_{\xi_1}^{u^1}
\sqrt{ b_1(\eta)/\Psi(\eta)}
\,d\eta
\quad
\right]_{\xi_1}
$
and
$
\left[
\int_{0}^{u^1}
\sqrt{ b_1(\eta)/\Psi(\eta)}
\,d\eta
\quad
\right]_{0}
$
are isomorphic: moreover so are, for each $k$,
those of
$
\left[
\int_{\xi_k}^{u^k}
\sqrt{
f_k(\eta)
\,d\eta
}
\right]_{\xi_k}
$
and
$
\left[
\int_{0}^{u^k}
\sqrt{
f_k(\eta)
\,d\eta
}
\right]_{0}
$.
\QUAN
\end{lemma}
\begin{theorem}
\label{teoremaprincipale}
\label{bla1}
A direct manifold
$
{\cal U} 
$
with metric
$
\Lambda(u^1.....u^N)=
b_1(u^1)\,du^1\odot du^1+
\sum_{i=2}^N a_i(u^1)f_i(u^i)
\,du^i\odot du^i
$
is geodesically complete if and only if it is coercive.
\end{theorem}
{\bf Proof:} 
a) suppose that ${\cal U}$ is {\TTT coercive} and that $U$ 
is an element of geodesic, starting at
a metrically ordinary point; moreover, let 
$
(\U^{\bullet}{}^1(0)...\U^{\bullet}{}^N(0)    )
$
be the initial velocity of $U$.
Suppose at first that $ z\mapsto u^1(z)$ is a constant function (hence $\U^{\bullet}{}^1(0)=0\,$): then, by lemma \ref{integraleprimo2}, the equations of $U$ are
\begin{equation}
\label{riportointpr2}
\cases{
{u}^1(z)
=A_1\quad\cr
(\U^{\bullet}{}^k(z)   )^2
f_k(u^k(z)   )
=A_k\quad\  k=2...N.
}
\end{equation}
The Riemann surface of
$ z\mapsto {u}^1(z)$ is trivially isomorphic to $ \PI^1$; if  $ A_k=0$ so is the one of $ z\mapsto {u}^k(z)$ is isomorphic to $ (\PI^1,\id,\id,A)$ for some complex constant $A$; if $ A_k\not=0$ we could rewrite the k-th equation of (\ref{riportointpr2}) in the form:
\begin{equation}
\label{riscritta}
\frac{1}{B_k}\int_{u^k(0)}^{u^{k}(z)} \phi(\eta)\,d\eta
=\,z,
\end{equation}
where $\phi_k^2=f_k$ and $B_k^2=A_k$, the choice of $\phi_k$ and $B_k$ being made in such a way that 
$ \U^{\bullet}{}^k(0){\phi_k(0)}={B_k}$.
By hypothesis, the Riemann surface $(S_k,\pi_k,j_k,\Phi_k   )$
of the HFG $\left[\int_0^{u^k}\phi_k\,d\eta\right]_{0}$
is such that $\PI^1\setminus\Phi_1(S_1)$ is a finite set; by lemma \ref{isom} {\bf 1)} 
the Riemann surface of the HFG $\left[\int_{u^k(0)}^{u^k}\phi_k\,d\eta\right]_{u^k(0)}$
is isomorphic to $(S_k,\pi_k,j_k,\Phi_k   )$; but, by (\ref{riscritta}), the germs
$ \u^k_{z=0}$ and
$\left[\int_{u^k(0)}^{u^k}\phi_k\,d\eta\right]_{u^k(0)}$
are each one inverse of the other; hence, by lemma \ref{inverse} the Riemann surface of $ \u^k_{z=0}$ is complete; this eventually implies that the Riemann surface of the element
$
z\longmapsto(u^1(z)...u^N(z)    )
$
is complete too: this fact ends the proof of a) in the case that $u^1$ is a constant function.
Otherwise, by lemma \ref{integraleprimo1}, the equations of $U$ are
\begin{equation}
\label{riportointpr1}
\cases{
(\U^{\bullet}{}^1(z)   )^2
( b_1(u^1(z)   )   )
=A_1-\sum_{l=2}^N\frac{A_l}
{a_l(u^1(z)   )}
\quad\spadesuit\cr
(\U^{\bullet}{}^k(z)   )^2
f_k(u^k(z)   )
\left[a_k(u^1(z)   )\right]^2
=A_k\quad k=2...N\ \clubsuit.}
\end{equation}
for suitable complex constants $A_1...A_N.$
Consider now the germ $ z\mapsto u^1(z)$ in $z=0$:
rewrite the first equation of (\ref{riportointpr1}) in the form:
\begin{equation}
\label{riscritta2}
\int_{u^1(0)}^{u^{1}(z)}
\frac{ d\eta}
 {\alph(\eta)_{u^1(0)}}
=\,z,
\end{equation}
where 
$
(
\alph(\eta)_{u^1(0)}
)^2=
({A_1-\sum_{l=2}^N
A_l
/
a_l(\eta)})
/
{b_1(\eta)}
$ 
in a neighbourhood of $z=0$,
the choice of the square root  $\alph_k$ being 
made in such a way that
$
 \alph_{u^1(0)}(u^1(0))=1/\U^{\bullet}{}^1(0)
$. 
Denote now  by $\alph_{u=0}$ the HFG 
such that
$
(
\alph_{0}
)^2=
\left[
\frac
{1}
{ b_1}
(
{
A_1-\sum_{l=2}^N\frac
{A_l}
{a_l}}
)
\right]_{0}
$,
the choice of the 'square root' $\alph_0$ being arbitrary.
By hypothesis, the Riemann surface $(S_1,\pi_1,j_1,\Phi_1   )$
of the HFG $\left[\int_0^{u^1}1/\alph_0\right]_0$
is such that $\PI^1\setminus\Phi_1(S_1)$ is a finite set.
By lemma \ref{isom} the Riemann surfaces of $\left[\int_0^{u^1}1/\alph_0\right]_0$ and of $\left[\int_{u^1_0}^{u^1}1/\alph_0\right]_{u^1_0}$
are both isomorphic to $(S_1,\pi_1,j_1,\Phi_1   )$; but, by (\ref{riscritta}), the germs
$ \u^1_{z=0}$ and
$[\int_0^{u^1}1/\alph_0]_{u^1(0)}$ are each one inverse of the other; hence, by lemma \ref{inverse} the Riemann surface of $ \u^1_{z=0}$ is complete.
Let now $ 2\leq k\leq N$: if $ A_k=0$ the Riemann surface of $ z\mapsto {u}^k(z)$ is isomorphic to $ (\PI^1,\id,\id,A)$ for some complex constant $A$; if $ A_k\not=0$ we could rewrite the k-th equation of (\ref{riportointpr1}) in the form:
\begin{equation}
\label{riscritta3}
\int_{u^k(0)}^{u^{k}(z)} \phi(\eta)\,d\eta
=\,
\int_0^z\frac{B_k\,dz}{a_k(u^1(z)   )},
\end{equation}
where $\phi_k^2=f_k$ and $B_k^2=A_k$, the choice of $\phi_k$ and $B_k$ being made in such a way that 
$ \U^{\bullet}{}^k(0)\,\phi(
u^k(0)   )\,
a_k(u^1(z)   )={B_k}$.
Denote now  by $[\varphi_k]_{u^k=0}$ the HFG 
defined by setting
$[\varphi_k]_{u^k=0}^2=
\left[
f_k
\right]_{u^k=0}
$,
the choice of the "square root" $[\varphi_k]_{u^k=0}$ being arbitrary.
By hypothesis, the Riemann surface $(S_k,\pi_k,j_k,\Phi_k   )$
of the HFG $\left[\int_0^{u^k}\varphi_k\right]_{0}$
is such that $\PI^1\setminus\Phi_1(S_1)$ is a finite set; moreover, by lemma \ref{isom} the Riemann surfaces of the HFG $\left[\int_{u^k(0)}^{u^k}\phi_k\,d\eta\right]_{u^k(0)}$
is isomorphic to $(S_k,\pi_k,j_k,\Phi_k   )$; but, by (\ref{riscritta3}) the germs
$ \left[z\longrightarrow\u^k\right]_{z=0}$,
$\left[\int_{u^k(0)}^{u^k}\phi_k\,d\eta\right]_{u^k(0)}$ and
$
\left[z\longrightarrow\int_0^z\frac{B_k}{a_k(u^1(\zeta)   )}\,d\zeta\right]_{z=0}$
satisfy, in the above order, the hypotheses of lemma \ref{quasiinverse}  {\bf 2)};
moreover, the Riemann surface with L-singularities of $\left[\int_{u^k(0)}^{u^k}\phi_k\,d\eta\right]_{u^k(0)}$ is complete, since the one of 
$\left[\phi_k\right]_{u^k(0)}$ is complete without L-singularities.
Therefore the Riemann surface with L-singularities of $ \u^k_{z=0}$ is 
complete, hence so is the one of 
$
z\longmapsto(u^1(z)...u^N(z)    ),
$:
 this fact ends the proof of a).
Vice versa, suppose that 
$
{\cal U}
$
is not {coercive}: then
{\BBB either}
there exists a complex
n-tuple $ (A_1...A_N   )\in\CI^N$ such that 
$
b_1(x_0^1)\not=0$,
$
A_1-\sum_{l=2}^N\frac
{
A_l}
{
a_l(x_0^1)}\not=0
$
and for each one of the two HFG's 
${\alph_1}$
and
${\alph_2}$
such that 
$
({\alph_i})^2=
\left[
\frac{1}{b_1}
(
A_1-\sum_{l=2}^N\frac
{A_l}
{a_l}
)
\right]
_{0}
\quad i=1,2
$,
the Riemann surface $(S_1,\pi_1,j_1,\Phi_1   )$
of both the HFG's
$
\left[
\int_{x_0}^{u^1}
{
\frac{d\,\eta}
{\alph_i(\eta)}
}
\right]_{x_0^1}
\quad i=1,2;
$ 
is such that $\PI^1\setminus\Phi_1(S_1)$ is
an infinite set;
{\BBB or}
there exists $k$, $2\leq k\leq N$
such that, for each one of the two HFG's 
$\left[\phi_{k1}\right]_0$
and
$\left[\phi_{k2}\right]_0$
such that 
$
\left[\phi_{ki}\right]_0=
\left[
f_k
\right]_0,
\quad i=(1,2)
$,
the Riemann surface $(S_k,\pi_k,j_k,\Phi_k   )$
of both the HFG's
$\left[\int_0^{u^k}
\phi_ki(\eta)
\,d\eta\right]_{0}\ \ i=1,2$
is such that $\PI^1\setminus\Phi_1(S_1)$ is an infinite set.
In the first case the geodesic element
$
z\longmapsto U(z)=(u^1(z)...u^N(z)    )
$
starting from $0$ with velocity 
$(L_1...L_N    )$, such that
$ 
L_1^2=
\frac
{1}
{b_1(0)}
(
A_1-\sum_{l=2}^N\frac
{A_l}
{a_l(0)}
)
,\ L_k^2=\frac{A_k}{f_k(0)a_k(0)},\ 
k=2...N
$,
satisfies the equation
$
\int_0^{u^1(z)}
\frac
{ d\eta}
{ \alph_i(\eta)}
=z
$,
where $i=1$ or $i=2$;
by lemma \ref{inverse}, this fact implies that the Riemann surface of $\left[
 z\longmapsto
u^1(z)
\right]_{0}  $
is incomplete, hence  the same holds about $ z\longmapsto
U(z)$.
Consider now the second case: first construct a geodesic element 
$
z\longmapsto U(z)=(0...u^k(z)...0    )
$
with all constant components except $ u^k, k\geq 2$.
Now recall lemma \ref{integraleprimo2} to conclude that 
$
z\longmapsto u^k(z)
$
satisfies, in a neighbourhood of $z=0$ the equation 
$
\frac
{1}
{C_k}
{\int_0^{u^k(z)} \phi_{ki}(\eta)\,d\eta}=
$,
for a suitable complex constant $A_k$; therefore its Riemann surface is incomplete by lemma \ref{inverse}; this fact ends the proof.
\QUAN
\begin{definition}
\label{directbiholom}
\rm
Let ${\cal U}$ and ${\cal V}$ be direct manifolds   they are {\BBB directly biholomorphic} provided that they are biholomorphic under a direct product of biholomorphic functions between each  ${\cal U}_i$ and each ${\cal V}_i$.
\end{definition}
\begin{remark}
\label{biholom}
\rm
Definition \ref{coercive} is invariant by direct biholomorphism (see definition \ref{directbiholom}
) : in other words, if ${\cal U}$ and ${\cal V}$ are directly biholomorphic, then 
${\cal U}$ is coercive if and only ${\cal V}$ is too: this is a simple consequence of 'changing variable' in integrals in definition \ref{coercive}.
\end{remark}
Therefore, we could yield the following
\begin{definition}
\label{equicoercive}
\rm
An equivalence class $\left[{\cal U}    \right]$ of direct manifolds, consisting of mutually directly (see definition \ref{directbiholom}
) biholomorphic elements is {\BBB coercive}
provided that any one of its representatives is coercive. 
\end{definition}
Our goal is now to extend definitions \ref{coercive} and \ref{equicoercive} to warped products 
containg some $\PI^1$'s among their factors.
Keeping into account remark \ref{biholom}, consider a warped product 
$(\prod_{i=1}^N {\cal U}_i,\Lambda)$
of Riemann spheres, complex planes or one-dimensional unit balls, which we shall call direct manifold too;
let $ L\subset \{1...N\}$ be the set of indices such that $
{\cal U}_l\simeq \PI^1$ for each $l\in L$.
\begin{definition}
\rm
Let $ Y=(y^1...y^N )\in{\cal U}$: then $(Y,L    )$ is a 
{\BBB principal multipole} of ${\cal U}$ provided that 
$
b_1(y^1)=\infty$ and 
$
f_l(y^l)=\infty$ for each $l\in L\setminus \{1\}$;
\label{multipole}
A direct manifold 
$(\prod_{i=1}^N {\cal U}_i,\Lambda)$
of Riemann spheres, complex planes or one-dimensional unit balls with metric is
{\BBB partially projective} if some one of its factors is biholomorphic to the Riemann sphere
$\PI^1$;
\label{multicoer}
a partially direct manifold $ {\cal U}$ is
{\BBB coercive in opposition to the principal multipole } $(Y,L    )$ if, set 
$
{\cal W}_i=
\cases
{
{\cal U}_i & if $i\not\in L$\cr
{\cal U}_i\setminus \{y^i\} & if $i\in L$,
}
$
then $ \prod_{i=1}^N {\cal W}_i$ is coercive in the sense of definition \ref{equicoercive}, that is to say, belongs to a coercive equivalence class with respect to direct biholomorphicity.
\end{definition}
\subsection{Warped product of Riemann surfaces}
Consider now the warped product of Riemann surfaces
$$
{\cal S}={\cal S}_1\times_{a_2}
{\cal S}_2\times_{a_3}
{\cal S}_3\times
........
\times_{a_N}{\cal S}_N,
$$
where each ${\cal S}_i$ is endowed with meromorphic metric $\lambda_i$: ${\cal S}$' metric $\Lambda$ is defined by setting
$
\Lambda=\lambda_1+\sum_{k=2}^N a_k\lambda_k
$,
and each $a_k$ is a not everywhere vanishing meromorphic function on ${\cal S}_i$:
as a simple consequence of Riemann's uniformization theorem,
${\cal S}$ admits universal covering $ \Psi : {\cal U}\longrightarrow {\cal S}$, where ${\cal U}$
is a direct
manifold, endowed with the pull-back meromorphic metric $ \Psi^*\Lambda$:
this universal covering is unique up to direct biholomorphisms.
\begin{definition}
\rm
\label{rieco}
${\cal S}$ is
{\BBB totally unelliptic} provided that none of the ${\cal S}_i$ is elliptic;
{$L$-elliptic} provided that there exists a nonempty 
set of indices $L$ such that ${\cal S}_l$ is elliptic if and only if $l\in L$.

If ${\cal S}$ is a $L$-elliptic warped product with universal covering 
$\Psi : {\cal U}\longrightarrow{\cal S}$, then $(Z,L    )$ is a {\BBB principal multipole }
for ${\cal S}$ provided that $Z\in {\cal S}$ and each $Y\in\Psi^{-1}(Z    )$ is a 
principal multipole for ${\cal U}$.

A totally unelliptic warped product of Riemann sur\-faces is {\BBB coercive} provided 
that its universal covering is coercive in the sense of definition \ref{equicoercive};
a $L$-elliptic warped product of Riemann surfaces is {\BBB coercive in opposition to the 
principal multipole} $(Z,L    )$ provided that its universal
covering ${\cal U}$ is coercive in opposition to each principal multipole
$(Y,L    )$ as $Y$ runs over $\Psi^{-1}(Z)$.
\end{definition}
\begin{theorem}
\label{unell}
\label{bla2}
A totally unelliptic warped product of Riemann surfaces ${\cal S}$ is geodesically complete if and only if it is coercive.
\end{theorem}
{\bf Proof:} let $\Psi : {\cal U}\longrightarrow
{\cal S}$ be the universal covering of ${\cal S}$:
by definition \ref{rieco} ${\cal U}$ is coercive, hence geodesically complete by theorem \ref{teoremaprincipale}.
Let now $\gam$ be a germ of geodesic in ${\cal S}$, starting at a metrically ordinary point: since $\Psi$ is a local isometry, there exists a germ $\bet$ of geodesic in ${\cal U}$, starting at a metrically ordinary point, such that $\gam=\Psi\circ\bet$.
By definition of completeness, the Riemann surface with L-singularities $(\Sigma,\pi,j,B,{\cal U}    )$ of $\bet$ is such that $\PI^1\setminus \pi(\Sigma    )$ is a finite set; moreover, $(\Sigma,\pi,j,\Psi\circ B,{\cal S}    )$ is an analytical continuation, with L-singularities, of $\gam$.
This proves that, if $(\widetilde \Sigma,\widetilde \pi,\widetilde j,G,{\cal S}    )$ is the Riemann surface with L-singularities of $\gam$, then $\PI^1\setminus \widetilde \pi(\widetilde \Sigma    )$ is a finite set too, hence ${\cal S}$ is geodesically complete.
On the other side, if ${\cal S}$ admits an incomplete germ of geodesic $\gam$, starting at a metrically ordinary point, then there exists  an incomplete germ of geodesic $\bet$ in ${\cal U}$, starting at a metrically ordinary point, such that 
$\gam=\Psi\circ\bet$; this means by theorem \ref{teoremaprincipale},
that ${\cal U}$ is not coercive; eventually, by definition \ref{rieco}, ${\cal S}$ is not coercive: this fact ends the proof.
\QUAN
\begin{theorem}
\label{bla3}
A $L$-elliptic warped product of Riemann
surfaces ${\cal S}$ is geodesically complete if
and only if it is coercive in opposition to some principal  multipole.
\end{theorem}
{\bf Proof:} suppose that ${\cal S}$ is coercive in opposition to some principal  multipole $(Z,L    )$: then, by theorem \ref{unell},  ${\cal S}$ is coercive in opposition to $(Z,L    )$ if and only if ${\cal S}\setminus Z$ is geodesically complete;
since $Z$ is not metrically ordinary, ${\cal S}$ is geodesically complete.
On the other hand, suppose that ${\cal S}$ admits an incomplete geodesic
$(\Sigma,\pi,j,\gamma,{\cal S}    )$: let $(Z,L    )$ be a principal multipole of ${\cal S}$ wich is known to exist; set $R=\gamma^{-1}({\cal S}\setminus Z    )\subset\Sigma$.
Now $(R,\pi\vert_R,j,\gamma\vert_R,{\cal S}\setminus Z    )$ is an incomplete geodesic of ${\cal S}\setminus Z$: this fact implies that ${\cal S}\setminus Z$ is not geodesically complete, hence it is not coercive, that is to say, ${\cal S}$ is not coercive in opposition to $(Z,L    )$.
The arbitrariness of $Z$ allows us to conclude the proof.
\QUAN
\subsection{Examples}
We show a wide class of coercive
direct manifolds.
To do this, we need some technicalities from integral calculus, hence we state:
\begin{proposition}
\label{integrale}\rm
Set $\Delta=b^2-4ac$, the germ $\left[
\frac{1}
{\sqrt{a\eta^2+b\eta+c}}
\right]_0 
$ admits one of the following primitives, depending on $a,b,c$:

$
\left[
\frac{1}{\sqrt{a}}\log(\eta+\frac{b}{2a}+
\sqrt{\eta^2+\frac{b}{a}\eta+\frac{c}{a}}    )+cost
\right]_0$
{the same branch of $\sqrt{\ }$, any branch of the logarithm},
{if $a\not=0$ and $\Delta\not=0$};

$
\left[
\frac{ 2}{ b}\sqrt{b\eta+c}+cost 
\right]_0
$
{the same branch of $\sqrt{\ }$},
{if $a=0$ and $b\not=0$};

$
\left[ 
\eta/\sqrt{c}+cost\right]_0
$
the same branch of $\sqrt{\ }$,
{if $a=b=0$.}
\end{proposition}
Let now $S_i,\ i=1..N$ be Riemann surfaces, which we suppose for simplicity 
parabolic or hyperbolic, $p_i\colon{\cal U}_i\longrightarrow S_i$ 
their universal covering, where each
${\cal U}_i\simeq\CI$ or $\DI$;
finally, let $\phi_i$ be meromorphic 
functions such that $\phi_1\circ p_1$ and $(\phi_i\circ p_i)^{\prime},\ i=1..N$ take all complex values but at most a finite number (the hypothesis on $phi_i\circ p_i$ could be weakened;  even dropped, if $S_i$ is parabolic: see \cite{hayman}, introduction).
Moreover, let $(a_i,\ b_i,\ c_i)\in\CI^3\setminus
0\ i=1..N $, set
$
S=\prod_{i=1}^N S_i,\ {\cal U}=\prod_{i=1}^N={\cal U}_i,\ p=(p_1....p_N) 
$ 
and consider the meromorphic metric
$$
\Lambda=d\phi_1\odot d\phi_1+\sum
_{i=1}^N \frac {d\phi_i\odot d\phi_i}{a_i\phi_1^2+b_i\phi_1+c_i}.
$$
\begin{theorem}
$({\cal U},\Lambda)$ is coercive (hence geodesically complete).
\end{theorem}
{\bf Proof}: by pulling back $\Lambda$ with respect to the universal covering $p$ we get
$$
p^*\Lambda(z^1...z^N)=
\left[(\phi_1\circ p_1)^{\prime}\right]^2dz^1\odot dz^1\ +\sum
_{i=1}^N \frac 
{\left[(\phi_i\circ p_i)^{\prime}\right]^2dz^i\odot dz^i}
{a_i(\phi_1\circ p_1)^2+b_i\phi_1\circ p_1+c_i}.
$$
We claim that $\left({\cal U},p^*\Lambda\right)$ is coercive: 
indeed, for every n-tuple $\displaystyle \left(A_1...A_N   \right)\in\CI^N$ such that 
$
{
(\phi_1\circ p_1)^{\prime}(0)\not=0}
$
and
$
{A_1-\sum_{l=2}^N
{A_l}
{a_i(\phi_1\circ p_1)^2+b_i\phi_1\circ p_1+c_i}\not=0,
}
$,
set $\phi\circ p_1=\psi$,
there holds
$$
\int_0^{u^1}
{
{\sssqrt{\displaystyle A_1-\sum_{l=2}^N
{A_l}
{(a_i(\psi)^2+b_i\psi+c_i)(\eta)}}}
}
{(\psi)^{\prime}(\eta)d\,\eta}
=\Phi\left(\psi\right)
,
$$
where $\Phi$ is one (depending on the constants $A_1...A_N$) of the holomorphic function 
germs on the right hand member of proposition \ref{integrale}.
 
This fact shows that the maximal analytical continuation of $u^1\longrightarrow\Phi
\left(\phi_1\circ p_1(u^1)    \right)$ takes all $\PI^1$'s values but a finite number, because so 
does the meromorphic function $\phi_1$ and hence $\phi_1\circ p_1$;
moreover,
for each $i$, $2\leq i\leq N$,
each one of the two HFG's 
$
\pm
\left[(\phi_i\circ p_i)^{\prime}\right]
$
could be continuated to $\pm \left[(\phi_i\circ p_i)^{\prime}\right]$ which, by 
assumption, takes all values but at most two ones.
\QUAN
%

%
\noindent
\scriptsize
In this paper we investigate possible extensions of the idea of
geodesic completeness in complex manifolds, following two
directions: metrics are somewhere allowed not to be of maximum
rank, or to have 'poles' somewhere else. Geodesics are
eventually defined on Riemann surfaces over regions in the
Riemann sphere.
Completeness theorems are given in the framework of warped
products of Riemann surfaces.
\end{document}